\documentclass[11pt]{article} 
\usepackage[utf8]{inputenc}
\usepackage{caption,graphicx,amsmath,amsthm,amssymb, authblk,xcolor,enumitem,subfigure,float,verbatim,hyperref,pgfplots,bibentry,mathtools,soul,isomath, cancel}
\usepackage{mathrsfs}
\usepackage{indentfirst}
\usepackage{nccmath}
\graphicspath{ {images/} }
\usepackage{multicol, MnSymbol}
\usepackage{graphicx,amsmath,amssymb,amsthm,marvosym,tikz,fontenc}
\usetikzlibrary{backgrounds}
\usetikzlibrary{patterns, shapes.misc, positioning}
\usepackage{pgfplots}
\pgfplotsset{width=10cm,compat=1.9,tick scale binop=\times}
\usepackage{indentfirst}
\usepackage[a4paper,bindingoffset=0.2in,%
left=0.8in,right=1in,top=1.2in,bottom=1.2in,%
footskip=.25in]{geometry}
\usepackage{relsize}
\usepackage[english]{babel}
\allowdisplaybreaks
\theoremstyle{plain}
\newtheorem{theorem}{Theorem}[section]
\newtheorem{lemma}[theorem]{Lemma}

\newtheorem{corollary}[theorem]{Corollary}
\newtheorem{definition}[theorem]{Definition}

\date{}
\begin{document}
\title
{\bf{Eccentricity spectral properties of $\mathcal{C}$-graphs}}
\author {\small Anjitha Ashokan
\footnote{anjitha\_p220118ma@nitc.ac.in}   and Chithra A V
\footnote{chithra@nitc.ac.in} \\ \small Department of Mathematics, National Institute of Technology Calicut,\\
\small Kerala, India-673601}
\date{}
\maketitle
\begin{abstract} 
A cograph is a simple graph that contains no induced path on four vertices. In this paper, we consider $\mathcal{C}$-graphs, which are a specific class of cographs, defined as
 $$\overline{\overline{\overline{K_{\alpha_{1}}}\cup K_{\alpha_{2}}}\cup \cdots \cup K_{\alpha_{2k}}},$$ 
 where $k \geq 2$, $\alpha_{2k} \geq 2$, and $K_{\alpha_{i}}$ denotes the complete graph on $\alpha_{i}$ vertices.
We investigate the spectral properties of the eccentricity matrix of this particular class of cographs.
Additionally, we determine the irreducibility and inertia of the eccentricity matrix of  $\mathcal{C}$-graphs. Furthermore, we identify an interval $(-1-\sqrt{2},-2)\cup (-2,0)$ in which these graphs have no eccentricity eigenvalues.
\\
   
\end{abstract}

\section{Introduction}

Let $G=(V(G),E(G))$ be a finite simple connected graph with vertex set $V(G)=\{v_{1},v_{2},\ldots, v_{n}\}$ and edge set $E(G)=\{e_{1},e_{2}.\ldots , e_{m}\}.$ If two vertices $v_{i}$ and $v_{j}$ are adjacent, this is denoted by $v_{i}\sim v_{j}.$
The \textit{adjacency matrix} $A(G)$ of $G$ is a $n\times n$ matrix whose rows and columns  are indexed by the vertex set of the graph and the entries are defined by \\
\begin{equation*}
   A(G)=(a_{ij})=\begin{cases}
            1 & \text{ if } v_{i}\sim v_{j}\\
            0 &\text{ otherwise.}
            \end{cases}
\end{equation*}
The \textit{distance} between the vertices $v_{i},v_{j}\in V(G),$ denoted by $d(v_{i},v_{j}),$ is defined to be the smallest value among the lengths of the paths between the vertices $v_{i}$ and $v_{j}.$ The \textit{distance matrix} $D(G)$ of  a connected graph $G$ is an $n\times n$ matrix whose $(i,j)th$ entry is equal to $d(v_{i},v_{j}).$ 
The \textit{eccentricity} $e(v_{i})$ of a vertex $v_{i}$ in a graph $G$ is defined as $e(v_{i})=max\{d(v_{i},v_{j}):v_{j}\in V(G)\}.$
The \textit{eccentricity matrix} $\epsilon(G)$ of a graph $G$ is obtained from its distance matrix $D(G)$ by retaining the largest distance in each row and each column, and setting the other entries in that row and column to zero. More precisely, the $(i,j)$th entry of the eccentricity matrix is defined as  \begin{equation*}
  ( \epsilon(G))_{ij}=\begin{cases}
            d(v_{i},v_{j}) & \text{ if } d(v_{i},v_{j})=\min\{e(v_{i}),e(v_{j})\},\\
            0 &\text{ otherwise.}
            \end{cases}
\end{equation*}
Randic$^{'}$ et al.\cite{MR3136762} defined the \textit{$D_{MAX}$ matrix}, which was renamed as the eccentricity matrix by Wang et al.\cite{MR3906706}.  Since $\epsilon(G)$ is a symmetric matrix, all its eigenvalues are real.
Let $\lambda_{1}(\epsilon(G)), \lambda_{2}(\epsilon(G)), \ldots, \lambda_{n}(\epsilon(G))$ be the eigenvalues of $\epsilon(G) $ (eccentricity eigenvalues) listed in non decreasing order, that is, $\lambda_{1}(\epsilon(G))\leq \lambda_{2}(\epsilon(G))\leq \cdots \leq \lambda_{n}(\epsilon(G)).$ 
The multiset of distinct eccentricity eigenvalues of $G$ is the eccentricity spectrum of $G,$  denoted by $Spec_{\epsilon}(G)=\{\lambda_{1}(\epsilon(G))^{m_{1}},\lambda_{2}(\epsilon(G))^{m_{2}}, \cdots \, \lambda_{s}(\epsilon(G))^{m_{s}}\},$ where 
$m_{i}$ is the multiplicity of $\lambda_{i}(\epsilon(G)).$
An $n\times n$  nonnegative matrix $N$ is \textit{reducible} if there exists an $n\times n$ permutation matrix $P$ such that $PNP^{T}=\begin{pmatrix}
    N_{11} & N_{12}\\
    0 & N_{22}
\end{pmatrix},$ where $N_{11}$ is a $r\times r$ submatrix with $1\leq r < n.$ If no such permutation matrix $P$ exists, then $N$ is \textit{irreducible}.
Unlike the distance and adjacency matrices, the eccentricity matrix is not always irreducible.
 
 Recently, the study of eccentricity matrices of graphs has attracted significant attention within the research community. 
 Several spectral properties of the eccentricity matrix have been explored (see \cite{MR4110109, MR4574900}  and references therein). In particular, recent investigations have focused on its spectral characteristics, with emerging applications in the field of chemistry \cite{MR3136762, MR3906706}.

A cograph is a graph with no induced $P_{4}.$  The diameter of a connected cograph is at most two.
Cographs were introduced in the 1960s \cite{kelmans1965number}, and many equivalent definitions for cographs are given in \cite{MR3442533}.
The structural and spectral properties, as well as various applications of cographs, have been studied  \cite{gagneur2004modular, MR3804786}.

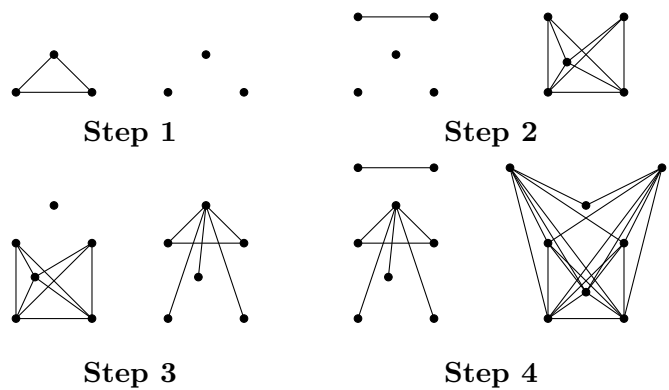
\begin{figure}[H]
 \centering
 \begin{tikzpicture}[scale=.5]
 \filldraw[fill=black](-20,0)circle(0.1cm); 
 \filldraw[fill=black](-22,0)circle(0.1cm); 
 \filldraw[fill=black](-21,1)circle(0.1cm); 
 \draw(-20, 0)--(-22,0);
 \draw(-20,0)--(-21,1);
  \draw(-22,0)--(-21,1);
 \filldraw[fill=black](-18,0)circle(0.1cm); 
 \filldraw[fill=black](-16,0)circle(0.1cm); 
 \filldraw[fill=black](-17,1)circle(0.1cm); 
 \node at (-19, -1.1) {\bfseries Step 1};
 \filldraw[fill=black](-11,0)circle(0.1cm); 
 \filldraw[fill=black](-13,0)circle(0.1cm); 
 \filldraw[fill=black](-12,1)circle(0.1cm); 
 \filldraw[fill=black](-13,2)circle(0.1cm);
 \filldraw[fill=black](-11,2)circle(0.1cm);
  \draw(-13,2)--(-11,2);
\filldraw[fill=black](-8,0)circle(0.1cm); 
 \filldraw[fill=black](-6,0)circle(0.1cm); 
 \filldraw[fill=black](-7.5,0.8)circle(0.1cm); 
 \filldraw[fill=black](-8,2)circle(0.1cm);
 \filldraw[fill=black](-6,2)circle(0.1cm);
 \draw(-8,0)--(-6,0);
 \draw(-8,0)--(-7.5,0.8);
\draw(-6,0)--(-7.5,0.8);
\draw(-8,2)--(-7.5,0.8);
\draw(-8,2)--(-6,0);
\draw(-8,2)--(-8,0);
\draw(-6,2)--(-7.5,0.8);
\draw(-6,2)--(-6,0);
\draw(-6,2)--(-8,0);
 \node at (-9.5, -1.1) {\bfseries Step 2};
\filldraw[fill=black](-22,-6)circle(0.1cm); 
 \filldraw[fill=black](-20,-6)circle(0.1cm); 
 \filldraw[fill=black](-21.5,-4.9)circle(0.1cm); 
 \filldraw[fill=black](-22,-4)circle(0.1cm);
 \filldraw[fill=black](-20,-4)circle(0.1cm);
\filldraw[fill=black](-21,-3)circle(0.1cm);
 \draw(-22,-6)--(-20,-6);
  \draw(-22,-6)--(-21.5,-4.9);
   \draw(-21.5,-4.9)--(-20,-6);
    \draw(-22,-4)--(-21.5,-4.9);
     \draw(-22,-4)--(-22,-6);
      \draw(-22,-4)--(-20,-6);
       \draw(-20,-4)--(-21.5,-4.9);
     \draw(-20,-4)--(-22,-6);
      \draw(-20,-4)--(-20,-6);
  \filldraw[fill=black](-18,-6)circle(0.1cm); 
 \filldraw[fill=black](-16,-6)circle(0.1cm); 
 \filldraw[fill=black](-17.2,-4.9)circle(0.1cm); 
 \filldraw[fill=black](-18,-4)circle(0.1cm);
 \filldraw[fill=black](-16,-4)circle(0.1cm);
\filldraw[fill=black](-17,-3)circle(0.1cm); 
 \draw(-17,-3)--(-18,-6);
  \draw(-17,-3)--(-16,-6);
   \draw(-17,-3)--(-17.2,-4.9);
    \draw(-17,-3)--(-18,-4);
     \draw(-17,-3)--(-16,-4);
     \draw(-18,-4)--(-16,-4);
      \node at (-19, -7.5) {\bfseries Step 3};
  \filldraw[fill=black](-13,-6)circle(0.1cm); 
 \filldraw[fill=black](-11,-6)circle(0.1cm); 
 \filldraw[fill=black](-12.2,-4.9)circle(0.1cm); 
 \filldraw[fill=black](-13,-4)circle(0.1cm);
 \filldraw[fill=black](-11,-4)circle(0.1cm);
\filldraw[fill=black](-12,-3)circle(0.1cm); 
\filldraw[fill=black](-13,-2)circle(0.1cm); 
\filldraw[fill=black](-11,-2)circle(0.1cm); 
 \draw(-12,-3)--(-13,-6);
  \draw(-12,-3)--(-11,-6);
   \draw(-12,-3)--(-12.2,-4.9);
    \draw(-12,-3)--(-13,-4);
     \draw(-12,-3)--(-11,-4);
     \draw(-13,-4)--(-11,-4);
   \draw(-13,-2)--(-11,-2);  
 \filldraw[fill=black](-8,-6)circle(0.1cm); 
 \filldraw[fill=black](-6,-6)circle(0.1cm); 
 \filldraw[fill=black](-7,-5.3)circle(0.1cm); 
 \filldraw[fill=black](-8,-4)circle(0.1cm);
 \filldraw[fill=black](-6,-4)circle(0.1cm);
\filldraw[fill=black](-7,-3)circle(0.1cm); 
\filldraw[fill=black](-9,-2)circle(0.1cm); 
\filldraw[fill=black](-5,-2)circle(0.1cm); 
\draw(-8,-6)--(-6,-6); 
 \draw(-7,-5.3)--(-8,-6); 
 \draw(-7,-5.3)--(-6,-6); 
  \draw(-8,-4)--(-8,-6);   
  \draw(-8,-4)--(-6,-6); 
  \draw(-8,-4)--(-7,-5.3); 
\draw(-6,-4)--(-8,-6);   
  \draw(-6,-4)--(-6,-6); 
  \draw(-6,-4)--(-7,-5.3); 
\draw(-9,-2)--(-8,-6); 
\draw(-9,-2)--(-6,-6); 
\draw(-9,-2)--(-8,-4); 
\draw(-9,-2)--(-6,-4); 
\draw(-9,-2)--(-7,-5.3); 
\draw(-9,-2)--(-7,-3); 
\draw(-5,-2)--(-8,-6); 
\draw(-5,-2)--(-6,-6); 
\draw(-5,-2)--(-8,-4); 
\draw(-5,-2)--(-6,-4); 
\draw(-5,-2)--(-7,-5.3); 
\draw(-5,-2)--(-7,-3); 
\node at (-9.5, -7.5) {\bfseries Step 4};
 \end{tikzpicture}
  \caption{$C(3,2,1,2)$}
 \label{fig:1}
 \end{figure}

 In 2024, Santanu et al.\cite{mandal2024laplacian} defined a particular type of cograph as follows. 
 Let $K_{\alpha}$ denote the complete graph on $\alpha$ vertices. For a sequence of  positive integers $\alpha_{1},\alpha_{2},\ldots,\alpha_{l},$  the cograph $C(\alpha_{1},\alpha_{2},\ldots,\alpha_{l})$ is defined recursively by \begin{align*}
            C(\alpha_{1})&=\overline{K_{\alpha_{1}}}\\
            C(\alpha_{1},\alpha_{2},\ldots,\alpha_{i})&=\overline{C(\alpha_{1},\alpha_{2},\ldots,\alpha_{i-1})\cup K_{\alpha_{i}}}, \text{ for } 2\leq i\leq l.
\end{align*} 
The tuple $(\alpha_{1}, \alpha_{2}, \ldots, \alpha_{l})$ is known as the generating sequence of the associated cograph. 
In other words, to construct the cograph $C(\alpha_{1},\alpha_{2},\ldots \alpha_{l}),$ follow this process:
\begin{itemize}
\item Begin with the complete graph $K_{\alpha_{1}}.$
\item At each subsequent step $ i$, $2\leq i\leq l$, take the disjoint union of  $K_{\alpha_{i}}$ with the graph obtained in the previous step,
 and then take the complement of that union.
 \item After processing all elements of the sequence in this manner, the resulting graph is $C(\alpha_{1},\alpha_{2},\ldots \alpha_{l}).$
\end{itemize}
That is,  $$C(\alpha_{1},\alpha_{2},\ldots,\alpha_{l})  \cong \overline{\overline{\overline{K_{\alpha_{1}}}\cup K_{\alpha_{2}}}\cup \cdots \cup K_{\alpha_{l}}}.$$ 
The number of vertices in  $C(\alpha_{1},\alpha_{2},\ldots \alpha_{l})$ is 
$\sum_{i=1}^{l}\alpha_{i}.$\\
\newpage
An illustration of a 8-vertex $\mathcal{C}$-graph $C(3,2,1,2)$ is shown in Figure.\ref{fig:1}.

Let $\mathcal{C}$ denote the class of cographs constructed as described above. If $G \in \mathcal{C}$, then $G$ is a $\mathcal{C}$-graph.
The $\mathcal{C}$-graphs $C(1,1,1,1,1)$ and $C(1,2,1,1)$ 
represent the same $\mathcal{C}$-graph. If $G\in \mathcal{C}$ is constructed in an even number of steps, then its representation is unique \cite{MR4854054}. 
This paper focuses on the $\mathcal{C}$-graphs where $l$ is even.

 A cograph is a threshold graph if it is a $\{P_{4}, 2K_{2}, C_{4}\}-$ free graph, that is, a particular cograph without an induced subgraph isomorphic to two parallel edges or the $4-$ vertex cycle or a path on $4-$ vertices.
Every threshold graph can be constructed from a finite binary sequence  \cite{MR4039442, MR4541509}. 
 In 2011, Sciriha and Farrugia \cite{sciriha2011spectrum} provided some spectral properties of the adjacency eigenvalues of threshold graphs, and some additional properties of the adjacency (resp. distance, eccentricity) eigenvalues of threshold graphs were explored in \cite{MR4541509, MR3116409, MR3809377, MR4364833}. 
Chain graphs are another class of graphs constructed from a finite sequence. In \cite{MR4338338}, the authors determined an eigenvalue-free interval of the distance matrix of chain graphs, and in \cite{MR4847141}, the eccentricity inertia of chain graphs has been computed.
In \cite{mandal2024laplacian, MR4854054}, spectral properties of $\mathcal{C}$-graphs related to the adjacency and Laplacian matrices were investigated.

Inspired by these works, in this paper, we study the eccentricity spectral properties  of $\mathcal{C}$-graph $G$ of the form
\begin{equation*}
    G=C(\alpha_{1},\alpha_{2},\ldots,\alpha_{2k}),\text{ with } k\geq 2,  \alpha_{2k}\geq 2.
\end{equation*}

This article is organized as follows. In Section $2,$ we recall some basic results along with some notations to make the paper more self-contained. In Section $3,$ we estimate the inertia  
for the eccentricity matrix of the $\mathcal{C}$-graphs. Also, we check the irreducibility of the eccentricity matrix of $\mathcal{C}$-graphs. In Section $4,$ we obtain an eccentricity eigenvalue free interval for $\mathcal{C}-$graphs.

\section{Preliminaries}
\begin{lemma}\cite{horn1994topics}\label{Weyls inequality}
Let $B,$ $C$ be Hermitian matrices of order $n$ and let the respective eigenvalues of $B, C,B+C$ be $\lambda_{i}(B),$ $\lambda_{i}(C),$ $\lambda_{i}(B+C),i=1,2,\ldots n$ each  ordered as: \begin{equation}\label{ordering of eigenvalues}
    \lambda_{1}\leq \lambda_{2}\leq \cdots \leq \lambda_{n}.
\end{equation}
 Then, 
 \begin{equation*}
    \lambda_{i}(B)+\lambda_{1}(C)\leq \lambda_{i}(B+C)\leq \lambda_{i}(B)+\lambda_{n}(C),
\end{equation*} 
where $i=1,2,\ldots, n.$
\end{lemma}



\begin{lemma}\cite{horn1994topics}\label{Sylverster's law of inertia}
    Let $B$ and $C$ be Hermitian matrices of equal order. Then $B$ and $C$ have the same inertia if and only if $B=SCS^{*}$ for some non-singulr matrix $S,$ where $S^{*}$ is the conjugate transpose of $S.$
\end{lemma}
\begin{lemma}\label{ostrowski theorem}\cite{horn1994topics}
    Let $B,$ $S$ be matrices of order $n$ with $B$ Hermitian and $S$ nonsingular. Let the eigenvalues of $B,$ $SBS^{*}$ and $SS^{*}$ be arranged in non-decreasing order. Let $\sigma_{1}\geq \sigma_{2}\geq \cdots \geq \sigma_{n}$ be the singular values of $S.$ Then for each $p=1,2,\ldots, n$ there is a positive real number $\theta_{p}\in [\sigma_{n}^{2}, \sigma_{1}^{2}]$ such that \begin{align*}
        \lambda_{p}(SBS^{*})=\theta_{p}\lambda_{p}(B).
    \end{align*} 
\end{lemma}

Let $N$ be a square symmetric matrix.
The inertia of a matrix $N$ is represented by a triple $(n_{-}(N), n_{0}(N), n_{+}(N)),$ where $n_{-}(N)$ (resp., $n_{0}(N), n_{+}(N))$ is the count of negative (resp., zero, positive) eigenvalues of $N.$
\\
An antiregular graph is a connected graph in which at most two vertices have the same degree. The antiregular graph on $n$ vertices is denoted by $H_{n}$. 
%
%
An antiregular graph is a $\mathcal{C}$-graph   and it's representation is either $C(1,1,1,\ldots,1)$ or $C(1,2,1,\ldots,1)$ \cite{MR4854054}.\\
Throughout this paper we denote the adjacency matrix of a $n-$ vertex antiregular graph by $A_{n}.$

\begin{lemma}\label{inertia of An}\cite{MR3116409}
    Let $H_{2k-1}$ be an antiregular graph on $2k-1$ vertices. Then the inertia of $A_{2k-1}$ is,
     $$ inertia (A_{2k-1})=(n_{-}(A_{2k-1}), n_{0}(A_{2k-1}), n_{+}(A_{2k-1}))=(\frac{2k-2}{2},1,\frac{2k-2}{2}).$$
\end{lemma}

For the matrix $A_{n}$, the eigenvalues $0$ and $-1$ are called the trivial eigenvalues, and all other eigenvalues are called non-trivial eigenvalues. 

 Let $\lambda_{+}(A_{n})$ denote the smallest positive eigenvalue of the matrix $A_{n},$ and let $\lambda^{-}(A_{n})$ denote the largest negative eigenvalue of $A_{n}$ smaller than $-1.$
\begin{lemma}\cite{MR4039442}\label{max and min eigenvalue bound for Hn}
    Let $n$ be a positive integer. Then 
      $$ \lambda^{-}(A_{n})<\frac{-1-\sqrt{2}}{2} \text{ and  } \lambda_{+}(A_{n})>\frac{-1+\sqrt{2}}{2}.$$
      Equivalently, the closed interval $[\frac{-1-\sqrt{2}}{2}, \frac{-1+\sqrt{2}}{2}]$ does not contain any non-trivial eigenvalue of the antiregular graph $H_{n}.$
\end{lemma}
\begin{lemma}\cite{MR2499317}\label{-1 is not an eigenvalue of A2k-1}
 For $n\geq 0,$ the value   $-1$ is not an eigenvalue for $A_{2n+1}$.
\end{lemma}

\begin{definition}(Equitable partition)\cite{MR1829620}
 Suppose that $M$ is an $n\times n$ real symmetric matrix whose rows and columns are indexed by $X=\{1,2,\ldots, n\}.$ Let $\Pi$ be a partition of $X,$  $X=X_{1}\cup X_{2} \cup \cdots \cup X_{m}$. The matrix $M$ can be written as 
\begin{equation*}
M=\begin{pmatrix}
    M_{1,1} & M_{1,2}&\cdots &M_{1,m}\\
    M_{2,1} &M_{2,2} &\cdots &M_{2,m}\\
    \vdots&\vdots&\ddots&\vdots\\
    M_{m,1}&M_{m,2}&\cdots&M_{m,m}
\end{pmatrix}
\end{equation*}
where $M_{i,j}$ is the submatrix of $M$ whose rows and columns are indexed by $X_{i}$ and $X_{j},$ respectively, for $1\leq i,j\leq m.$ Let $b_{ij}$ be the average row sum of $M_{i,j}.$ Then $Q_{m}=(b_{ij})_{m\times m}$ is the quotient matrix of $M$ with respect to the partition $\Pi.$ If the row sum of each block $M_{i,j}$ is a constant, then the partition is an equitable partition. 
\end{definition}
Note that the quotient matrix $Q_{m}$ may not be symmetric even for equitable partitions.
\begin{lemma}\cite{MR1829620}\label{Q as equitable quotient matrix}
Let $M$ be a real symmetric matrix and $\Pi$ be an equitable partition of $M$ with quotient matrix $Q_{m}$. Then the characteristic polynomial of the quotient matrix $Q_{m}$ divides the characteristic polynomial of $M.$
    
\end{lemma}

\begin{lemma}\label{interlacing thm}(Interlacing theorem)\cite{horn1994topics}
    Let $M$ be a Hermitian matrix of order $n,$ and let $H$ be a principal submatrix of $M$ of order $m.$ If $\theta_{1}(M)\leq \theta_{2}(M)\leq \cdots \leq \theta_{n}(M)$ lists the eigenvalues of $M$ and $\mu_{1}(H)\leq\mu_{2}(H)\leq\cdots \leq \mu_{m}(H)$ the eigenvalues of $H,$ then $\theta_{i}(M)\leq \mu_{i}(H)\leq \theta_{i+n-m}(M)$ for $1\leq i\leq m.$ 
\end{lemma}
The principal submatrix of a matrix 
$M$ is obtained by
deleting the same set of rows and the corresponding columns from $M.$
\begin{lemma}\label{result for the existence of principle submatrix}\cite{MR4777066}
 Let $G$ be a connected graph and $H$ is a subgraph of $G.$  If $e_{H}(v_{i})=e_{G}(v_{i})$ holds for every $v_{i}\in V(H)$ and $d_{G}(v_{i},v_{j})=d_{H}(v_{i},v_{j})$ holds for every pair $v_{i},v_{j}\in V(H),$ then $\epsilon(H)$ is a principal submatrix of $\epsilon(G).$
\end{lemma}
 
For any connected graph $G,$ the \textit{eccentric graph}, $G^{e},$ is defined as the graph with $V(G^{e})=V(G)$ and any two vertices $u, v\in V(G^{e})$ are adjacent if and only if the distance detween them is min$\{e(u),e(v)\}$ in $G.$
\begin{lemma}\cite{MR4522999}\label{for irreducibility}
    Let $G$ be a  $(p,q)$ graph. Then the matrix $\epsilon(G)$ is irreducible if and only if $G^{e}$ is a connected graph. 
\end{lemma}

\section{The eccentricity eigenvalues of \texorpdfstring{$\mathcal{C}$}{C}-graphs}

Let $G=C(\alpha_{1},\alpha_{2},\ldots,\alpha_{2k})$.
%
The vertex set of $G$  can be partitioned as   $V(G)=V_{1}\cup V_{2} \cup \cdots \cup  V_{2k},$  where $V_{i}$ denotes the set of $\alpha_{i}$ vertices introduced in the $i-$th step.
Note that $G$ has a diameter of $2.$ \\
The following theorem determines whether the eccentricity matrix of a $\mathcal{C}$-graph is irreducible or reducible, depending entirely on the final parameter $\alpha_{2k}$.
\begin{theorem}
    Let $G=C(\alpha_{1},\alpha_{2},\ldots,\alpha_{2k}), k\geq 2$. Then
    \begin{enumerate} 
        \item If $\alpha_{2k}=1$, then $\epsilon(G)$ is irreducible 
        \item If $\alpha_{2k} \geq 2$, then $\epsilon(G)$ is reducible. 
    \end{enumerate}
    \end{theorem}
\begin{proof}
    Let $V(G)=V_{1} \cup V_{2} \cup \cdots \cup V_{2k},$ be the vertex set partition of $G,$ also $|V_{i}|=\alpha_{i}, 1\leq i\leq 2k.$ 
    \begin{enumerate}
        \item If $\alpha_{2k}=1,$ and $u$ be the vertex in $V_{2k}.$ Then, $e(u)=1.$ For every  $ v\in V(G)-\{V_{2k}\},$ we have
        $d(u,v)=1=min\{e(u),e(v)\} .$ Hence,  $G^{e}$ is a connected graph.By Lemma \ref{for irreducibility}, $\epsilon(G)$ is irreducible.
        \item If $\alpha_{2k} \geq 2,$ then $e(v)=2$ for every $v \in V_{2k}.$ For any $w \in V(G)-V_{2k},$ we have $d(w,v)=1 \neq min\{e(w),e(v)\}.$ Hence, no edge  exists between $V_{2k}$ and $V_{i}, \text{ for } 1\leq i \leq 2k-1.$ Therefore, $G^{e}$ is disconnected.  By Lemma \ref{for irreducibility}, $\epsilon(G)$ is reducible.
    \end{enumerate}
\end{proof}


For $k=1$, the $\mathcal{C}$-graph is given by $G=C(\alpha_{1},\alpha_{2})=K_{\alpha_{1}} \vee \alpha_{2}K_{1}$,  where $\vee$ denotes the join of the graphs $K_{\alpha_{1}}$ and $\alpha_{2}K_{1}$. 
Hence, we have the following result.
\begin{lemma}\cite{MR2523015}
   If $G=K_{\alpha_{1}} \vee \alpha_{2}K_{1},$ then  $$Spec_{\epsilon}(G)=\begin{pmatrix}
                                                                           -1&-2&\alpha_{1}+\alpha_{2}-2-\frac{\alpha_{1}-1}{2}\pm \sqrt{(\alpha_{1}-\alpha_{2}-\frac{\alpha_{1}-1}{2})^{2}+\alpha_{1}\alpha_{2}}\\
                                                                           \alpha_{1}-1&\alpha_{2}-1&1
                                                                         \end{pmatrix}.$$
\end{lemma}

\vspace{0.7cm}
This study focuses on the  eccentricity spectra  of $\mathcal{C}$-graphs of the form,
$$G=C(\alpha_{1},\alpha_{2},\ldots,\alpha_{2k}),\text{ with }   \alpha_{2k}\geq 2,k\geq 2.$$

The eccentricity matrix of $G$ can be written as

$\epsilon(G)=\scriptsize{\begin{pmatrix}
0_{\alpha_{1}\times\alpha_{1}} & 0_{\alpha_{1}\times\alpha_{2}} & 2J_{\alpha_{1}\times\alpha_{3}} & 0_{\alpha_{1}\times\alpha_{4}} & \cdots & 2J_{\alpha_{1}\times\alpha_{2k-1}} & 0_{\alpha_{1}\times\alpha_{2k}}\\
0_{\alpha_{2}\times\alpha_{1}} & 2(J-I)_{\alpha_{2}\times\alpha_{2}} & 2J_{\alpha_{2}\times\alpha_{3}} & 0_{\alpha_{2}\times\alpha_{4}} &  \cdots & 2J_{\alpha_{2}\times\alpha_{2k-1}} & 0_{\alpha_{1}\times\alpha_{2k}}\\
2J_{\alpha_{3}\times\alpha_{1}} & 2J_{\alpha_{3}\times\alpha_{2}} & 0_{\alpha_{3}\times\alpha_{3}} & 0_{\alpha_{3}\times\alpha_{4}} &  \cdots & 2J_{\alpha_{3}\times\alpha_{2k-1}} & 0_{\alpha_{3}\times\alpha_{2k}}\\
0_{\alpha_{4}\times \alpha_{1}}&0_{\alpha_{4}\times \alpha_{2}}&0_{\alpha_{4}\times \alpha_{3}}&2(J-I)_{\alpha_{4}\times \alpha_{4}}&\cdots &2J_{\alpha_{4}\times \alpha_{2k-1}}&0_{\alpha_{4}\times \alpha_{2k}}\\
 \vdots&\vdots&\vdots&\vdots&\ddots&\vdots&\vdots\\
 2J_{\alpha_{2k-1}\times \alpha_{1}}&2J_{\alpha_{2k-1}\times \alpha_{2}}&2J_{\alpha_{2k-1}\times \alpha_{3}}&2J_{\alpha_{2k-1}\times \alpha_{4}}&\cdots&0_{\alpha_{2k-1}\times \alpha_{2k-1}}&0_{\alpha_{2k-1}\times \alpha_{2k}}\\
 0_{\alpha_{2k}\times \alpha_{1}}&0_{\alpha_{2k}\times \alpha_{2}}&0_{\alpha_{2k}\times \alpha_{3}}&0_{\alpha_{2k}\times \alpha_{4}}&\cdots&0_{\alpha_{2k}\times \alpha_{2k-1}}&2(J-I)_{\alpha_{2k}\times \alpha_{2k}}
\end{pmatrix}}.$
\vspace{0.5cm}\\
Clearly, $\epsilon(G)$ is a symmetric matrix of order $n=\sum_{i=1}^{2k}\alpha_{i}.$

\begin{theorem}\label{multiplicity of 0 and -2}
    Let $G=C(\alpha_{1}, \alpha_{2},\ldots, \alpha_{2k}),$ $\alpha_{2k}\geq 2,$ $k\geq 2.$ Then
    \begin{enumerate}
        \item $-2$ is an eigenvalue of $\epsilon(G)$ with multiplicity at least $\alpha_{2}+\alpha_{4}+\cdots +\alpha_{2k}-k.$
        \item  $0$ is an eigenvalue of $\epsilon(G)$ with multiplicity at least $\alpha_{1}+\alpha_{3}+\cdots +\alpha_{2k-1}-k.$
         \item $2(\alpha_{2k}-1)$ is an eigenvalue of $\epsilon(G)$ with multiplicity at least $1.$
    \end{enumerate}
\end{theorem}
\begin{proof}
    Let $V(G)=V_{1} \cup V_{2}\cup \cdots\cup V_{2k-1} \cup V_{2k},$ where $V_{j}=\{v_{j,1}, v_{j,2},\ldots v_{j,\alpha_{j}}\}, \text{ where }1\leq j \leq 2k.$\\

For $i=1,2,\ldots k,$ define a vector $x^{(2i,l)} \in R^{n}$ to be 
\begin{equation*}
\begin{split}
 x^{(2i,l)}=(x_{v_{1,1}},x_{v_{1,2}},\ldots,x_{v_{1,\alpha_{1}}},x_{v_{2,1}},x_{v_{2,2}},\ldots,x_{v_{2,\alpha_{2}}},x_{v_{3,1}},x_{v_{3,2}},\ldots,x_{v_{3,\alpha_{3}}},&\\ x_{v_{4,1}},x_{v_{4,2}},\ldots,x_{v_{4,\alpha_{4}}},\cdots, x_{v_{2i,1}},x_{v_{2i,2}},\ldots,x_{v_{2i,\alpha_{2i}}},\cdots,x_{v_{2k,1}},x_{v_{2k,2}},\ldots,x_{v_{2k,\alpha_{2k}}})
 \end{split}
\end{equation*}
such that $x_{v_{2i,1}}=1,$ $x_{v_{2i,l}}=-1,$ for $ l=2,3,\ldots, \alpha_{2i}$ and $x_{v_{s,t}}=0,$ for $v_{s,t}\neq v_{2i,1},v_{2i,l}.$
Then, we have $\epsilon(G) x^{(2i,l)}=-2 x^{(2i,l)},$ and $ x^{(2i,l)}$ are linearly independent for $1\leq i\leq k$ and $2\leq l \leq \alpha_{2i}.$
 Thus, $-2$ is an eigenvalue of $\epsilon(G)$ with multiplicity at least $\alpha_{2}+\alpha_{4}+\cdots +\alpha_{2k}-k.$ \\
 
\vspace{0.5cm}
For $j=1,2,\ldots k,$ define a vector $y^{(2j-1,s)} \in R^{n}$  to be 
\begin{equation*}
\begin{split}
 y^{(2j-1,s)}=(y_{v_{1,1}},y_{v_{1,2}},\ldots,y_{v_{1,\alpha_{1}}},y_{v_{2,1}},y_{v_{2,2}},\ldots,y_{v_{2,\alpha_{2}}},y_{v_{3,1}},y_{v_{3,2}},\ldots,y_{v_{3,\alpha_{3}}},&\\ y_{v_{4,1}},y_{v_{4,2}},\ldots,y_{v_{4,\alpha_{4}}},\cdots, y_{v_{2j-1,1}},y_{v_{2j-1,2}},\ldots,y_{v_{2j-1,\alpha_{2j-1}}},\cdots,y_{v_{2k,1}},y_{v_{2k,2}},\ldots,y_{v_{2k,\alpha_{2k}}})
 \end{split}
\end{equation*}
such that $y_{v_{2j-1,1}}=1,$ $y_{v_{2j-1,s}}=-1,$ for $ s=2,3,\ldots, \alpha_{2j-1}$ and $y_{v_{p,q}}=0,$ for $v_{p,q}\neq v_{2j-1,1},v_{2j-1,s}.$
Then, we have $\epsilon(G) y^{(2j-1,s)}=0 y^{(2j-1,s)},$ and $ y^{(2j-1,s)}$ are linearly independent for $1\leq j\leq k$ and $2\leq s \leq \alpha_{2j-1}.$
 Thus, $0$ is an eigenvalue of $\epsilon(G)$ with multiplicity at least $\alpha_{1}+\alpha_{3}+\cdots +\alpha_{2k-1}-k.$
\\


Let $z $ be a vector in $R^{n}$ defined as,
$z=\begin{pmatrix}
    \mathbf{0}_{\alpha_{1}\times 1 }&\mathbf{0}_{\alpha_{2}\times 1}&\mathbf{0}_{\alpha_{3}\times 1}\cdots &\mathbf{0}_{\alpha_{2k-1}\times 1}&\mathbf{1}_{\alpha_{2k}\times 1}
\end{pmatrix}.$
Then, we have $\epsilon(G)z=2(\alpha_{2k}-1)z.$ Hence $2(\alpha_{2k}-1)$ is an eigenvalue of $\epsilon(G)$ with multiplicity atleast 1. 
\end{proof}

Next, we determine the other eccentricity eigenvalues of $G$ by using the equitable partition of $\epsilon(G)$. Consider the vertex partition $\Pi$ of $V,$ $V=V_{1}\cup V_{2}\cup \cdots \cup V_{2k-1} \cup V_{2k},$ where each $V_{i}$ denotes the set of $\alpha_{i}$ vertices introduced in the $i-$th step during the construction of  $G$. That is, 
$$V_{i}=\{v_{i,1}, v_{i,2},\ldots v_{i,\alpha_{i}}\},1\leq i \leq 2k. $$
Obviously, $\Pi$ is an equitable partition of $\epsilon(G),$ with quotient matrix $Q_{2k}(\epsilon(G)).$ 
\begin{align}\label{Q(G)}
       Q_{2k}(\epsilon(G))=&\begin{pmatrix}
   0&0&2\alpha_{3}&0&2\alpha_{5}&0&\cdots&2\alpha_{2k-1}&0\\
   0&2(\alpha_{2}-1)&2\alpha_{3}&0&2\alpha_{5}&0&\cdots&2\alpha_{2k-1}&0\\
   2\alpha_{1}&2\alpha_{2}&0&0&2\alpha_{5}&0&\cdots&2\alpha_{2k-1}&0\\
   0&0&0&2(\alpha_{4}-1)&2\alpha_{5}&0&\cdots&2\alpha_{2k-1}&0\\
   2\alpha_{1}&2\alpha_{2}&2\alpha_{3}&2\alpha_{4}&0&0&\cdots&2\alpha_{2k-1}&0\\
   0&0&0&0&0&2(\alpha_{6}-1)&\cdots&2\alpha_{2k-1}&0\\
   \vdots&\vdots&\vdots&\vdots&\vdots& \vdots&\ddots&\vdots&\vdots\\
2\alpha_{1}&2\alpha_{2}&2\alpha_{3}&2\alpha_{4}&2\alpha_{5}&2\alpha_{6}&\cdots&0&0\\
   0& 0& 0& 0& 0& 0&\cdots&0&2(\alpha_{2k}-1)
\end{pmatrix}.
\end{align}
Note that $Q_{2k}(\epsilon(G))$ is a matrix of order $2k,$ which is not necessarily symmetric.
By Lemma \ref{Q as equitable quotient matrix}, the eigenvalues of $Q_{2k}(\epsilon(G))$ are also eigenvalues of the eccentricity matrix   $\epsilon(G).$ 
 Moreover, from the structure of  $Q_{2k}(\epsilon(G))$, it is evident that
\begin{equation}
Spec(Q_{2k}(\epsilon(G)))=Spec(\tilde{Q}(\epsilon(G))) \cup \{2(\alpha_{2k}-1)\} \label{spectrum of Q(G) as union of Q tilde and one more },
\end{equation}
where $\tilde{Q}(\epsilon(G))$ denotes the submatrix of $Q_{2k}(\epsilon(G))$ obtained by removing the last row and last column of $Q_{2k}(\epsilon(G))$.
Let $D$ be a diagonal matrix defined as, $D=diag \begin{pmatrix}
\alpha_{1}&\alpha_{2}&\ldots& \alpha_{2k-1}
\end{pmatrix},$
 and a matrix $R$ is  defined as, 
\begin{align}\label{introducing R}
R=&D^{\frac{1}{2}}\tilde{Q}(\epsilon(G))D^{-\frac{1}{2}}\\
 =&\begin{pmatrix}
        0&0&2\sqrt{\alpha_{1}\alpha_{3}}&0&\cdots &2\sqrt{\alpha_{1}\alpha_{2k-1}}\\
        0&2(\alpha_{2}-1)&2\sqrt{\alpha_{2}\alpha_{3}}&0&\cdots& 2\sqrt{\alpha_{2}\alpha_{2k-1}}\\
        2\sqrt{\alpha_{1}\alpha_{3}}&2\sqrt{\alpha_{2}\alpha_{3}}&0&0&\cdots& 2\sqrt{\alpha_{3}\alpha_{2k-1}}\\
        0&0&0&2(\alpha_{4}-1)&\cdots& 2\sqrt{\alpha_{4}\alpha_{2k-1}}\\
\vdots&\vdots&\vdots&\vdots&\ddots&\vdots\\
2\sqrt{\alpha_{1}\alpha_{2k-1}}&2\sqrt{\alpha_{2}\alpha_{2k-1}}&2\sqrt{\alpha_{3}\alpha_{2k-1}}&2\sqrt{\alpha_{4}\alpha_{2k-1}}&\cdots &0\\
 \end{pmatrix}. 
\end{align}
Then, 
\begin{equation}\label{similarity between R and tilde Q}
D^{-\frac{1}{2}}RD^{\frac{1}{2}}=\tilde{Q}(\epsilon(G)).
\end{equation}
Hence, $R$ and $\tilde{Q}(\epsilon(G))$ are similar matrices.
Since $R$ is a symmetric matrix, $\tilde{Q}(\epsilon(G))$ is also a symmetric matrix. Thus, $\tilde{Q}(\epsilon(G))$ is a diagonalizable matrix.\\
 Next, we analyze the  spectral properties of  $\tilde{Q}(\epsilon(G))$ to obtain the main results of this section.

\begin{lemma}\label{multiplicity of the eigenvalue 0}
    Let $G=C(\alpha_{1},\alpha_{2},\ldots, \alpha_{2k}), \alpha_{2k}\geq 2$ $    k\geq 2$ be a $\mathcal{C}$-graph of order $n$. Then, $0$ is a simple eigenvalue of $Q_{2k}(\epsilon(G))$ if and only if $\alpha_{2}=1.$ 
\end{lemma}
\begin{proof}
To prove this, it suffices to prove that the rank of $\tilde{Q}(\epsilon(G))$ is less than $2k-1$ if and only if $\alpha_{2}=1.$
    Consider $\tilde{Q}(\epsilon(G))$ and perform the following row operations,
    \begin{enumerate}    
        \item $R_{2m+1}\rightarrow R_{2m+1}-R_{2m-1}, m=1, 2, \ldots , k-1$\\
        $\tilde{Q}(\epsilon(G))\sim \begin{pmatrix}
            0&0&2\alpha_{3}&0&2\alpha_{5}&0&\cdots&0&2\alpha_{2k-1}\\
            0&2(\alpha_{2}-1)&2\alpha_{3}&0&2\alpha_{5}&0&\cdots&0&2\alpha_{2k-1}\\
            2\alpha_{1}&2\alpha_{2}&-2\alpha_{3}&0&0&0&\cdots &0&0\\
            0&0&0&2(\alpha_{4}-1)&2\alpha_{5}&0&\dots&0&2\alpha_{2k-1}\\
            0&0&2\alpha_{3}&2\alpha_{4}&-2\alpha_{5}&0&\cdots &0&0\\
            0&0&0&0&0&2(\alpha_{6}-1)&\cdots&0&2\alpha_{2k-1}\\
            \vdots& \ \vdots& \vdots& \vdots& \vdots& \vdots& \ddots& \vdots& \vdots\\
            0&0&0&0&0&0&\cdots&2(\alpha_{2k-2}-1)&2\alpha_{2k-1}\\
            0&0&0&0&0&0&\cdots&2\alpha_{2k-2}&-2\alpha_{2k-1}
        \end{pmatrix}$
        \item $R_{1}\rightarrow R_{1}-R_{2}$

$\begin{pmatrix}
            0&-2(\alpha_{2}-1)&0&0&0&0&\cdots&0&0\\
            0&2(\alpha_{2}-1)&2\alpha_{3}&0&2\alpha_{5}&0&\cdots&0&2\alpha_{2k-1}\\
            2\alpha_{1}&2\alpha_{2}&-2\alpha_{3}&0&0&0&\cdots &0&0\\
            0&0&0&2(\alpha_{4}-1)&2\alpha_{5}&0&\dots&0&2\alpha_{2k-1}\\
            0&0&2\alpha_{3}&2\alpha_{4}&-2\alpha_{5}&0&\cdots &0&0\\
            0&0&0&0&0&2(\alpha_{6}-1)&\cdots&0&2\alpha_{2k-1}\\
            \vdots& \ \vdots& \vdots& \vdots& \vdots& \vdots& \ddots& \vdots& \vdots\\
            0&0&0&0&0&0&\cdots&2(\alpha_{2k-2}-1)&2\alpha_{2k-1}\\
            0&0&0&0&0&0&\cdots&2\alpha_{2k-2}&-2\alpha_{2k-1}
        \end{pmatrix}$

        \item $R_{2m}\rightarrow R_{2m}-R_{2m+2}, m= 1,2,\ldots, k-2$\\
        $\begin{pmatrix}
           0&-2(\alpha_{2}-1)&0&0&0&0&\cdots&0&0\\
           0&2(\alpha_{2}-1)&2\alpha_{3}&-2(\alpha_{4}-1)&0&0&\cdots&0&0\\
           2\alpha_{1}&2\alpha_{2}&-2\alpha_{3}&0&0&0&\cdots&0&0\\
           0&0&0&2(\alpha_{4}-1)&2\alpha_{5}&-2(\alpha_{6}-1)&\cdots&0&0\\
           0&0&2\alpha_{3}&2\alpha_{4}&-2\alpha_{5}&0&\cdots&0&0\\
           \vdots& \ \vdots& \ \vdots& \vdots& \vdots& \vdots& \ddots& \vdots& \vdots& \\
           0&0&0&0&0&0&\cdots&2(\alpha_{2k-2}-1)&2\alpha_{2k-1}\\
           0&0&0&0&0&0&\cdots&2\alpha_{2k-2}&-2\alpha_{2k-1}
        \end{pmatrix}$
        \item $R_{2}\leftrightarrow R_{3}$, $R_{4}\leftrightarrow R_{5}$, $R_{6}\leftrightarrow R_{7}, \ldots ,$ $R_{2k-2}\leftrightarrow R_{2k-1}.$ 
    \end{enumerate}
    After performing above operations we obtain a tridiagonal matrix, say $T,$ which is row equivalent to $\tilde{Q}(\epsilon(G)),$  \\
    $T=\begin{pmatrix}
        0&-2(\alpha_{2}-1)&0&0&0&0&\cdots&0&0&0\\
        2\alpha_{1}&2\alpha_{2}&-2\alpha_{3}&0&0&0&\cdots&0&0&0\\
        0& 2(\alpha_{2}-1)&2\alpha_{3}&-2(\alpha_{4}-1)&0&0&\cdots&0&0&0\\
        0 &0&2\alpha_{3}&2\alpha_{4}&-2\alpha_{5}&0&\cdots&0&0&0\\
        0&0&0&2(\alpha_{4}-1)&2\alpha_{5}&-2(\alpha_{6}-1)&\cdots&0&0&0\\
         \vdots&\vdots&\vdots&\vdots&\vdots& \vdots&\ddots&\vdots&\vdots\\
         0&0&0&0&0&0&\cdots&2\alpha_{2k-3}&2\alpha_{2k-2}&-2\alpha_{2k-1}\\
         0&0&0&0&0&0&\cdots&0&2(\alpha_{2k-2}-1)&2\alpha_{2k-1}
    \end{pmatrix}.$
    Thus, $rank(T)=2k-2$ if and only if $\alpha_{2}=1.$ 
    Hence, $0$ is a simple eigenvalue of $\tilde{Q}(\epsilon(G))$ if and only if $\alpha_{2}=1.$

    
\end{proof}

\begin{lemma}\label{multiplicity of eigenvalue -2}
    If $\lambda$ is an eigenvalue of $Q_{2k}(\epsilon(G)),$ then  $\lambda \neq -2.$
\end{lemma}
\begin{proof}
To prove this, it suffices to prove that the rank of $\tilde{Q}(\epsilon(G))+2I$ is equal to  $2k-1.$
    Consider $\tilde{Q}(\epsilon(G))+2I$ and perform the following row operations step by step,
    \begin{enumerate}    
        \item $R_{2m+1}\rightarrow R_{2m+1}-R_{2m-1}, m=1, 2, \ldots , k-1$ 
        \item $R_{1}\rightarrow R_{1}-R_{2}$
        \item $R_{2m} \rightarrow R_{2m}-R_{2m+2}, m= 1,2,\ldots, k-2$ 
        \item $R_{2}\leftrightarrow R_{3}$, $R_{4}\leftrightarrow R_{5}$, $R_{6}\leftrightarrow R_{7}, \ldots ,$ $R_{2k-2}\leftrightarrow R_{2k-1}.$
    \end{enumerate}
    After performing the above operations, we obtain a tridiagonal matrix, say $S,$ which is row equivalent to $\tilde{Q}(\epsilon(G))+2I,$ where 
    \begin{align*}
    S=\begin{pmatrix}
        2&-2\alpha_{2}&0&0&0&0&\cdots&0\\
        2(\alpha_{1}-1)&2\alpha_{2}&2(1-\alpha_{3})&0&0&0&\cdots&0\\
        0& 2\alpha_{2}&2\alpha_{3}&-2\alpha_{4}&0&0&\cdots&0\\
        0 &0&2(\alpha_{3}-1)&2\alpha_{4}&2(1-\alpha_{5})&0&\cdots&0\\
         \vdots&\vdots&\vdots&\vdots&\vdots& \vdots&\ddots&\vdots\\
         0&0&0&0&0&0&\cdots&2\alpha_{2k-1}
    \end{pmatrix}.
    \end{align*}
    
    Here, $S$ is a tridiagonal matrix. Moreover, $rank(S)=2k-1.$ Hence, $-2$ is not an  eigenvalue of $\tilde{Q}(\epsilon(G)).$ 
\end{proof}
\begin{lemma}\label{inertia of equitable quotient matrix}
Let $G=C(\alpha_{1},\alpha_{2},\ldots, \alpha_{2k})$, $k\geq 2,$ $ \alpha_{2k}\geq 2$ be a $\mathcal{C}$-graph. Then the inertia of the quotient matrix $Q_{2k}(\epsilon(G))$ is given by\\
\begin{enumerate}
    \item $(n_{-}(Q_{2k}(\epsilon(G))), n_{0}(Q_{2k}(\epsilon(G))), n_{+}(Q_{2k}(\epsilon(G)))=(k-1,0,k+1),$ if $\alpha_{2}\neq 1$
    \item $(n_{-}(Q_{2k}(\epsilon(G))), n_{0}(Q_{2k}(\epsilon(G))), n_{+}(Q_{2k}(\epsilon(G)))=(k-1,1,k),$ if $\alpha_{2}=1.$
\end{enumerate}
\end{lemma}
\begin{proof}
    Let $D=diag\begin{pmatrix}
    \alpha_{1}&\alpha_{2}&\ldots& \alpha_{2k-1} 
    \end{pmatrix}$ and 
    $\tilde{D}= diag \begin{pmatrix}
        0&\frac{\alpha_{2}-1}{\alpha_{2}}&0&\frac{\alpha_{4}-1}{\alpha_{4}}&\ldots&\frac{\alpha_{2k-2}-1}{\alpha_{2k-2}}&0
    \end{pmatrix},$  be two diagonal matrices of order $2k-1.$ \\
    Then, \begin{align}
         D^{\frac{1}{2}}2(A_{2k-1}+\tilde{D})D^\frac{1}{2}=&\begin{pmatrix}
            0&0&2\sqrt{\alpha_{1}\alpha_{3}}&0&\cdots &2\sqrt{\alpha_{1}\alpha_{2k-1}}\\
        0&2(\alpha_{2}-1)&2\sqrt{\alpha_{2}\alpha_{3}}&0&\cdots& 2\sqrt{\alpha_{2}\alpha_{2k-1}}\\
        2\sqrt{\alpha_{1}\alpha_{3}}&2\sqrt{\alpha_{2}\alpha_{3}}&0&0&\cdots& 2\sqrt{\alpha_{3}\alpha_{2k-1}}\\
        0&0&0&2(\alpha_{4}-1)&\cdots& 2\sqrt{\alpha_{4}\alpha_{2k-1}}\\
\vdots&\vdots&\vdots&\vdots&\ddots&\vdots\\
2\sqrt{\alpha_{1}\alpha_{2k-1}}&2\sqrt{\alpha_{2}\alpha_{2k-1}}&2\sqrt{\alpha_{3}\alpha_{2k-1}}&2\sqrt{\alpha_{4}\alpha_{2k-1}}&\cdots &0\\  \end{pmatrix}
         =&R,\label{equation connecting to antiregular graph}
        \end{align}
    where $A_{2k-1}$ denotes the adjacency matrix of the antiregular graph $C(1,1,\ldots, 1).$
    By $(\ref{similarity between R and tilde Q}),$  $\tilde{Q}(\epsilon(G))$ is similar to $D^{\frac{1}{2}}2(A_{2k-1}+\tilde{D})D^\frac{1}{2}.$\\
    By Lemma \ref{Sylverster's law of inertia}, $R$ and $2(A_{2k-1}+\tilde{D})$ have the same inertia. Hence,  $\tilde{Q}(\epsilon(G))$ and $2(A_{2k-1}+\tilde{D})$ have the same inertia.
 From Lemma \ref{inertia of An},
    it is known that inertia of $A_{2k-1}$ is $(n_{-}(A_{2k-1}), n_{0}(A_{2k-1}), n_{+}(A_{2k-1}))=(\frac{2k-2}{2},1,\frac{2k-2}{2}).$\\
 Now by Lemma \ref{Weyls inequality}, we have,
\begin{align*}
        \lambda_{k+1}(A_{2k-1}+\tilde{D}) & \geq \lambda_{k+1}(A_{2k-1})+\lambda_{1}(\tilde{D}) >0, \text{ and }\\
        \lambda_{k}(A_{2k-1}+\tilde{D})&\geq \lambda_{k}(A_{2k-1})+\lambda_{1}(\tilde{D}) =\lambda_{k}(A_{2k--1})=0,\\
        \lambda_{k}(A_{2k-1}+\tilde{D)}&\leq \lambda_{k}(A_{2k-1})+\lambda_{2k-1}(\tilde{D})=\lambda_{2k-1}(\tilde{D})<1.
        \end{align*}
 Hence,  $$0\leq \lambda_{k}(A_{2k-1}+\tilde{D})< 1.$$
        
  Now by Lemmas \ref{max and min eigenvalue bound for Hn} and  \ref{-1 is not an eigenvalue of A2k-1}, $\lambda_{k-1}(A_{2k-1}) < -1.$ Therefore, 
  \begin{align*}
      \lambda_{k-1}(A_{2k-1}+\tilde{D})&\leq \lambda_{k-1}(A_{2k-1})+\lambda_{2k-1}(\tilde{D}) \\
      &\leq -1+\lambda_{2k-1}(\tilde{D})\\
      &<0.
 \end{align*}

 \begin{enumerate}
     \item If $\alpha_{2}\neq 1,$ by Lemma \ref{multiplicity of the eigenvalue 0}, $0$ is not an eigenvalue of $\tilde{Q(\epsilon(G))},$ so is $A_{2k-1}+\tilde{D}.$\\
     Thus, $\lambda_{k}(A_{2k-1}+\tilde{D})>0.$
     Hence, inertia of $\tilde{Q}(\epsilon(G))
     =(k-1,0,k).$\\
     Therefore, inertia of $Q_{2k}(\epsilon(G))=
     (k-1,0,k+1).$
     \item If $\alpha_{2}=1,$ by Lemma \ref{multiplicity of the eigenvalue 0}, $0$ is  an eigenvalue of $\tilde{Q}(\epsilon(G)),$ so is $A_{2k-1}+\tilde{D}.$ \\
     Thus, $\lambda_{k}(A_{2k-1}+\tilde{D})=0.$
     Hence, inertia of $\tilde{Q}(\epsilon(G))=
     (k-1,1,k-1).$ \\
     Therefore, inertia of $Q_{2k}(\epsilon(G))=
     (k-1,1,k).$ 
 \end{enumerate}

\end{proof}

The following Theorem is a consequence of Theorem \ref{multiplicity of 0 and -2} and Lemmas \ref{multiplicity of the eigenvalue 0}, \ref{multiplicity of eigenvalue -2} and \ref{inertia of equitable quotient matrix}.
\begin{theorem}\label{inertia of C graph}
    Let $G=C(\alpha_{1},\alpha_{2},\ldots, \alpha_{2k})$, $k\geq 2,$ $ \alpha_{2k}\geq 2$ be a $\mathcal{C}$-graph. Then, the inertia of $\epsilon(G)$ is given by,
    \begin{enumerate}
        \item $(n_{-}(\epsilon(G)), n_{0}(\epsilon(G)), n_{+}(\epsilon(G))=(-1+\sum_{i=1}^{k}\alpha_{2i},-k+\sum_{i=1}^{k}\alpha_{2i-1}, k+1),$ if $\alpha_{2}\neq 1$
    \item $(n_{-}(\epsilon(G)), n_{0}(\epsilon(G)), n_{+}(\epsilon(G))=(-1+\sum_{i=1}^{k}\alpha_{2i}, -k+1+\sum_{i=1}^{k}\alpha_{2i-1}, k),$ if $\alpha_{2}=1.$
    \end{enumerate}
\end{theorem}
By applying Theorem \ref{multiplicity of 0 and -2}, Lemmas \ref{multiplicity of the eigenvalue 0} and \ref{multiplicity of eigenvalue -2}, we obtain the following result that specifies the exact multiplicities of the eigenvalues $0$ and $-2.$
\begin{corollary}\label{exact multiplicity of -2 and 0}
    Let $m_{0}(\epsilon(G))$ and  $m_{-2}(\epsilon(G))$ denote the multiplicities of the eigenvalues $0$ and $-2$, respectively, of the matrix $\epsilon(G).$ Then, 
    \begin{align*}
        m_{-2}(\epsilon(G))&= \sum_{i=1}^{k}\alpha_{2i}-k.\\
        m_{0}(\epsilon(G))&=\begin{cases}
            \sum_{i=1}^{k}\alpha_{2i-1}-k & \text{ ,if } \alpha_{2}=1 \\
            \sum_{i=1}^{k}\alpha_{2i-1}-k+1 &\text{ ,if } \alpha_{2}\neq 1.
            \end{cases}
\end{align*}
\end{corollary}
 The maximum number of distinct eigenvalues of the matrix $\epsilon(G)$ is determined by the eigenvalues $-2,$ $0,$ and those of the matrix $Q_{2k}(\epsilon(G)).$
 Hence, we obtain the following corollary.
\begin{corollary}
Let $m(\epsilon(G))$ denote the number of distinct eigenvalues of the matrix $\epsilon(G).$ Then, $m(\epsilon(G))\leq 2k+2.$
\end{corollary}

\section{Eccentricity eigenvalue-free interval for \texorpdfstring{$\mathcal{C}$}{C}-graphs}
In this section, we determine an interval in which the $\mathcal{C}$-graphs have no eccentricity eigenvalues.
\begin{lemma}\label{interlacing of Q}
    Let $C$ be a principal submatrix of $Q_{2k}(\epsilon(G))$ of size $t$.  
 The eigenvalues of $Q_{2k}(\epsilon(G))$ are denoted by $\lambda_1, \lambda_2, \dots, \lambda_{2k}$, ordered as 
$\lambda_1 \leq \lambda_2 \leq \cdots \leq \lambda_{2k}$, and the eigenvalues of $C$ by 
$\mu_1, \mu_2, \dots, \mu_t$, ordered as $\mu_1 \leq \mu_2 \leq \cdots \leq \mu_t$.  
Then, for $1 \leq i \leq t$, the following interlacing inequality holds: $\lambda_{i}\leq \mu_{i}\leq \lambda_{i+2k-t}$ for $1\leq i\leq m.$
\end{lemma}
\begin{proof}
 Let $D=diag(\alpha_{1},\alpha_{2},\ldots, \alpha_{2k})$ the diagonal matrix with $i-$ th diagonal entry $\alpha_{i}, 1\leq i\leq 2k.$
Then, we have \begin{align*}
    R=D^{\frac{1}{2}}Q_{2k}(\epsilon(G))D^{-\frac{1}{2}}&=\begin{pmatrix}
        0&0&2\sqrt{\alpha_{1}\alpha_{3}}&0&\cdots &2\sqrt{\alpha_{1}\alpha_{2k-1}}&0\\
        0&2(\alpha_{2}-1)&2\sqrt{\alpha_{2}\alpha_{3}}&0&\cdots& 2\sqrt{\alpha_{2}\alpha_{2k-1}}&0\\
        2\sqrt{\alpha_{1}\alpha_{3}}&2\sqrt{\alpha_{2}\alpha_{3}}&0&0&\cdots& 2\sqrt{\alpha_{3}\alpha_{2k-1}}&0\\
        0&0&0&2(\alpha_{4}-1)&\cdots& 2\sqrt{\alpha_{4}\alpha_{2k-1}}&0\\
\vdots&\vdots&\vdots&\vdots&\ddots&\vdots&\vdots\\
2\sqrt{\alpha_{1}\alpha_{2k-1}}&2\sqrt{\alpha_{2}\alpha_{2k-1}}&2\sqrt{\alpha_{3}\alpha_{2k-1}}&2\sqrt{\alpha_{4}\alpha_{2k-1}}&\cdots &0&0\\
0&0&0&0&\cdots &0&2(\alpha_{2k}-1)
          \end{pmatrix}.
\end{align*}
Therefore, $R$ and $Q_{2k}(\epsilon(G))$ are similar and have the same eigenvalues. 
Let $C$ be a principal submatrix of $Q_{2k}(\epsilon(G))$ of size $t.$ Without loss of generality, assume that $C$ is induced by the first $t$ rows and columns of $Q_{2k}(\epsilon(G)).$ Therefore we have $PQ_{2k}(\epsilon(G))P^{T}=C,$ where $P=(I_{t}|0)_{t\times 2k}$ and $I_{t}$ is the identity  matrix of order $t$.
Let \begin{align*}
    \tilde{C}&=PRP^{T}\\
        &=PD^{\frac{1}{2}}Q_{2k}(\epsilon(G))D^{-\frac{1}{2}}P^{T}\\
        &=PD^{\frac{1}{2}}P^{T}PQ_{2k}(\epsilon(G))P^{T}PD^{-\frac{1}{2}}P^{T}\\
        &=(PD^{\frac{1}{2}}P^{T})(PQ_{2k}(\epsilon(G))P^{T})(PD^{-\frac{1}{2}}P^{T})\\
        &=(PD^{\frac{1}{2}}P^{T})C(PD^{-\frac{1}{2}}P^{T}).
\end{align*}
Since $(PD^{-\frac{1}{2}}P^{T})=(PD^{\frac{1}{2}}P^{T})^{-1},$ we see that $\tilde{C}$ is similar to $C$, hence \textbf{$\mu_{1}\leq \mu_{2}\leq \cdots \leq \mu_{t}$ }are eigenvalues of $\tilde{C}.$
By Lemma \ref{interlacing thm} we have $\lambda_{i}\leq \mu_{i} \leq \lambda_{i+2k-t},$ $ 1\leq i\leq t.$
\end{proof}

We now investigate the important problem of identifying an eccentricity eigenvalue-free interval for  $\mathcal{C}$-graphs. Specifically, we aim to determine an interval that contains no eigenvalues of the matrix associated with $\mathcal {C}-$ graphs. Let $\lambda^{-}(\epsilon(G))$ denote  the largest negative eigenvalue of $\epsilon(G)$ smaller than $-2.$

\begin{theorem}\label{eigenvalue free interval for Q(G)}
  There is no eigenvalue of $Q_{2k}(\epsilon(G))$ lies in the interval $(-1-\sqrt{2},0).$  
\end{theorem}

\begin{proof}
 By  Lemma \ref{inertia of equitable quotient matrix}, it is enough to prove that $\lambda_{k-1}(Q_{2k}(\epsilon(G)))<-1-\sqrt{2}.$  Consider the $\mathcal{C}$-graph $G_{1}=C(1,1,\ldots, 1,2)$ on $2k+1$ vertices. Note that $G_{1}$ is a subgraph of $G.$  By Lemma \ref{result for the existence of principle submatrix}, $\epsilon(G_{1})$ is a principal submatrix of $\epsilon(G).$ Then,  Lemma \ref{interlacing thm},
\begin{equation}
                    \lambda_{k-1}(\epsilon(G))\leq \lambda_{k-1}(\epsilon(G_{1})).
   \label{equation from interlacing between G and G1} 
   \end{equation}
Now for the graph  $G_{1}$, using (\ref{introducing R}) and  (\ref{equation connecting to antiregular graph}), we have,
                    $$ \tilde{Q}(\epsilon(G_{1})) \sim 2A_{2k-1},$$
where $A_{2k-1}$ is the adjacency matrix of the antiregular graph $C(1,1,\ldots,1)$ on $2k-1$ vertices.\\
Now, \begin{align*}
    \lambda_{k-1}(\tilde{Q}(\epsilon(G_{1})))=&\lambda_{k-1}(2A_{2k-1})\\
                               =&2\lambda_{k-1}(A_{2k-1}).
       \end{align*}                        
      Using Lemmas \ref{inertia of An}, \ref{max and min eigenvalue bound for Hn} and \ref{-1 is not an eigenvalue of A2k-1} we have\\
      $$\lambda_{k-1}(A_{2k-1})<\frac{-1-\sqrt{2}}{2}.$$
      Hence,
      \begin{equation*}
                \lambda_{k-1}(\tilde{Q}(\epsilon(G_{1}))) < -1-\sqrt{2}.\\
\end{equation*}
 Since for any $\mathcal{C}$-graph  $G,$ the negative eigenvalues of $Q_{2k}(\epsilon(G))$ and $\tilde{Q}(\epsilon(G))$ are the same, it follows that, \\ 
 \begin{equation}\label{bound for Q(G1)}
            \lambda_{k-1}(Q_{2k}(\epsilon(G_{1}))) < -1-\sqrt{2}.
      \end{equation}  
     From Theorem \ref{inertia of C graph} and Lemma \ref{inertia of equitable quotient matrix}  we have, 
    \begin{equation*}
   (n_{-}(\epsilon(G_{1})), n_{0}(\epsilon(G_{1})), n_{+}(\epsilon(G_{1}))=(k,1,k), \\
   \end{equation*}
    and \\
   \begin{equation*}
   (n_{-}(Q_{2k}(\epsilon(G_{1}))), n_{0}(Q_{2k}(\epsilon(G_{1}))), n_{+}(Q_{2k}(\epsilon(G_{1})))=(k-1,1,k).
   \end{equation*}\\ 
   Also, $$Spec(\epsilon(G))=Spec(Q_{2k}(\epsilon(G)))\cup \{-2\}.$$
   Moreover,  $-2$ is not an eigenvalue of $Q_{2k}(\epsilon(G_{1})),$ and we have (\ref{bound for Q(G1)}). Hence,
   \
   $$\lambda_{k-1}(\epsilon(G_{1}))=\lambda_{k-1}(Q_{2k}(\epsilon(G_{1}))).$$
   Thus,  from (\ref{equation from interlacing between G and G1}), 
  \begin{equation} 
   \lambda_{k-1}(\epsilon(G))< -1-\sqrt{2}.\label{equation}
  \end{equation} 
    By Lemmas \ref{inertia of equitable quotient matrix}, and \ref{multiplicity of eigenvalue -2}, 
      $\lambda_{k-1}(\epsilon(G))=\lambda_{k-1}(Q_{2k}(\epsilon(G))).$
    Hence from (\ref{equation}), \\
  $$ \lambda_{k-1}(Q_{2k}(\epsilon(G)))<-1-\sqrt{2}.$$
\end{proof}


 Theorem \ref{eigenvalue free interval for Q(G)},  helps us to conclude that, $\lambda^{-}(\epsilon(G))=\lambda_{k-1}(Q_{2k}(\epsilon(G))).$ Now the following Corollary is a  direct consequence of Theorem \ref{eigenvalue free interval for Q(G)} and Theorem \ref{multiplicity of 0 and -2}.
\begin{corollary}\label{eigen value free interval of e(G)}
    Let $G=C(\alpha_{1},\alpha_{2},\ldots, \alpha_{2k}),k\geq 2,\alpha_{2k}\geq 2$ be a $\mathcal{C}$-graph then the interval $(-1-\sqrt{2}, -2)\cup (-2,0)$ does not contain any eigenvalue of $\epsilon(G).$
\end{corollary}
Combining Lemmas  \ref{multiplicity of the eigenvalue 0} and \ref{multiplicity of eigenvalue -2}, Corollary \ref{exact multiplicity of -2 and 0}, Theorem \ref{eigenvalue free interval for Q(G)}, and    Corollary \ref{eigen value free interval of e(G)}, we obtain the following result.
\begin{theorem}
Let $G=C(\alpha_{1},\alpha_{2},\ldots, \alpha_{2k}),k\geq 2,\alpha_{2k}\geq 2$    be a $\mathcal{C}$-graph.
\begin{enumerate}
 \item  If $\alpha_{2}\neq 1,$ then
 \begin{equation*}
     \begin{split}
          Spec_{\epsilon}(G)=&\{\lambda_{1}(Q_{2k}(\epsilon(G))), \lambda_{2}(Q_{2k}(\epsilon(G))),\ldots, \lambda_{k-1}(Q_{2k}(\epsilon(G))),-2^{\alpha_{2}+\alpha_{4}+\cdots+\alpha_{2k}-k},\\& 0^{\alpha_{1}+\alpha_{3}+\cdots+\alpha_{2k-1}-k},\lambda_{k}(Q_{2k}(\epsilon(G))),\lambda_{k+1}(Q_{2k}(\epsilon(G))),\ldots, \lambda_{2k}(Q_{2k}(\epsilon(G)))\},
\end{split}
          \end{equation*}
          where, $\lambda_{1}(Q_{2k}(\epsilon(G)))\leq \lambda_{2}(Q_{2k}(\epsilon(G)))\leq  \lambda_{k-1}(Q_{2k}(\epsilon(G)))<-1-\sqrt{2}<-2<0<\lambda_{k}(Q_{2k}(\epsilon(G))\leq \lambda_{k+1}Q_{2k}(\epsilon(G))\leq \cdots \leq \lambda_{2k}(Q_{2k}(\epsilon(G))). $
 \item If $\alpha_{2}=1,$ then 
    \begin{equation*}
    \begin{split}
    Spec_{\epsilon}(G)=&\{ \lambda_{1}(Q_{2k}(\epsilon(G))), \lambda_{2}(Q_{2k}(\epsilon(G))),\ldots, \lambda_{k-1}(Q_{2k}(\epsilon(G))),-2^{\alpha_{2}+\alpha_{4}+\cdots+\alpha_{2k}-k},\\& 0^{\alpha_{1}+\alpha_{3}+\cdots+\alpha_{2k-1}-k+1},\lambda_{k+1}(Q_{2k}(\epsilon(G))),\lambda_{k+2}(Q_{2k}(\epsilon(G))),\ldots, \lambda_{2k}(Q(\epsilon(G)))\},
\end{split}
    \end{equation*}
    where, $\lambda_{1}(Q_{2k}(G))\leq \lambda_{2}(Q_{2k}(G))\leq  \lambda_{k-1}(Q_{2k}(G))<-1-\sqrt{2}<-2<0<\lambda_{k+1}(Q_{2k}(G)\leq \lambda_{k+2}(Q_{2k}(G))\leq \cdots \leq \lambda_{2k}(Q_{2k}(G)). $
\end{enumerate}
\end{theorem}
Let $G=C(\alpha_{1},\alpha_{2},\ldots,\alpha_{2k}),\text{ with }   \alpha_{2k}\geq 2,k\geq 2$ be a $\mathcal{C}$-graph, with order $n$. Table \ref{table1} provides the eccentricity spectrum for some  $\mathcal{C}$-graphs.
\begin{table}[H]
\centering
\begin{tabular}{ccc}
\hline
 n & $G$  & $Spec_{\epsilon}(G)$ \\ 
\hline
5   & $\mathcal{C}(1,1,1,2)$   & $\{-2.8284,-2,0,2,2.8284\}$   \\ 
6   & $\mathcal{C}(2,1,1,2)$    & $\{-3.4641,-2,0^{2},2,3.4641\}$   \\ 
6  & $\mathcal{C}(1,1,2,2)$   & $\{-4,-2,0^{2},2,4\}$   \\
6 &  $\mathcal{C}(1,2,1,2)$   &  $\{-2.9624, -2,-2, 0.6222, 2, 4.3402\}$                       \\
7   & $\mathcal{C}(1,1,1,1,1,2)$  & $\{-3.4982,-2.5427,-2,0,0.6698,2,5.3711\}$   \\ 
7   & $\mathcal{C}(1,1,3,2)$   & $\{-4.8990,-2,0^{3},2,4.8990\}$  \\ 
7  & $\mathcal{C}(1, 3, 1, 2)$   & $\{-3.0283, -2^{3}, 0.8560, 2,6.123  \}$  \\ 
 7  & $\mathcal{C}(1,2,2,2)$   & $\{-4.35488,-2^{2},0,0.6433,2,5.7115 \}$   \\ 
7   & $\mathcal{C}(2,1,1,3)$   & $\{-3.4641,-2^{2},0^{2},3.4641,4\}$   \\ 
7  & $\mathcal{C}(1,2,1,3)$   & $\{-2.9624, -2^{3},0.6222, 4,4.3402 \}$  \\  
7  & $\mathcal{C}(1,1,1,4)$   &  $\{-2.8284, -2^{3}, 0, 2.8284,6 \}$   \\
\hline
\end{tabular}
\caption{The eccentricity spectrum of some $\mathcal{C}$-graphs.}
\label{table1}
\end{table}

\section*{Conclusion}
 In this paper, the eccentricity spectral properties of $\mathcal{C}$-graphs are studied. The irreducibility and inertia of the eccentricity matrix of  $\mathcal{C}$-graphs are also examined. Moreover, we showed that the interval $(-1 - \sqrt{2}, -2) \cup (-2, 0)$ contained no $\epsilon$-eigenvalues of the $\mathcal{C}$-graphs.

\section*{Acknowledgement}
The research of Anjitha Ashokan is supported by the University Grants Commission of India under the beneficiary code BININ05086971.
\section*{Declarations}
 On behalf of all authors, the corresponding author states that there is no conflict of interest.

\bibliographystyle{plain}
 \bibliography{reference}
\end{document}